\documentclass[]{interact}

\usepackage{pdflscape}

\usepackage{fullpage}
\usepackage{amsmath,amssymb,amsthm}
\usepackage{hyperref}
\usepackage{algorithm}
\usepackage{graphicx}
\usepackage{color}

\usepackage{pifont}
\newcommand{\cmark}{\ding{51}}%

\newcommand{\bs}{\boldsymbol}
\newcommand{\norm}[1]{\lVert #1 \rVert}
\DeclareMathOperator{\prox}{prox}
\DeclareMathOperator{\dom}{dom}
\DeclareMathOperator{\dist}{dist}
\DeclareMathOperator*{\locmin}{loc-min}

\newcommand{\revisionone}[1]{#1}

\title{A generic coordinate descent solver for nonsmooth convex optimization}

\author{\name{Olivier Fercoq\thanks{Email: olivier.fercoq@telecom-paristech.fr}}
\affil{LTCI, T\'el\'ecom ParisTech, Universit\'e Paris-Saclay, 46 rue Barrault, 75634 Paris Cedex 13, France}}

\begin{document}
\maketitle

\begin{abstract}
We present a generic coordinate descent solver for the minimization of a nonsmooth convex objective with structure. The method can deal in particular with problems with linear constraints. The implementation makes use of efficient residual updates and automatically determines which dual variables should be duplicated. A list of basic functional atoms is pre-compiled for efficiency and a modelling language in Python allows the user to combine them at run time. So, the algorithm can be used to solve a large variety of problems including Lasso, sparse multinomial logistic regression, linear and quadratic programs.
\end{abstract}


\section{Introduction}

Coordinate descent methods decompose a large optimization problem
into a sequence of one-dimensional optimization problems.
The algorithm was first described for the minimization of 
quadratic functions by Gauss and Seidel in~\cite{seidel1874ueber}.
Coordinate descent methods have become unavoidable in machine learning
because they are very efficient for key problems, namely Lasso~\cite{Friedman_Hastie_Hofling_Tibshirani07},
logistic regression~\cite{yu2011dual} and support vector machines~\cite{platt199912,SDCA}.
Moreover, the decomposition into small subproblems means that
only a small part of
the data is processed at each iteration and this makes coordinate descent
easily scalable to high dimensions.

One of the main ingredients of an efficient coordinate descent solver is 
its ability to compute efficiently partial derivatives of the objective function~\cite{nesterov2012efficiency}. In the case of least squares for instance, this  involves the definition of a vector of residuals that will be updated 
during the run of the algorithm. 
As this operation needs to be performed at each iteration, and millions of iterations are usually needed, the residual update and directional derivative computation must be coded in a compiled programming language.

Many coordinate descent solvers have been written in order to solve a
large variety of problems. However, most of the existing solvers can only solve problems of the type
\[
\min_{x \in \mathbb R^N} \sum_{j=1}^J f(A_j x - b_j) + \sum_{i=1}^I g(x^{(i)})
\]
where $x^{(i)}\in \mathbb R^{N_i}$ is the $i$th block of $x$, $\sum_{i=1}^I N_i = N$, $A_j \in \mathbb R^{M_j \times N}$ is a matrix and $b_j \in \mathbb R^{M_j}$ is a vector, and where $f$ is a convex differentiable function and $g$ is a convex lower-semicontinuous function whose proximal operator is easy to compute (a.k.a. a proximal-friendly convex function).
Each piece of code usually covers only one type of function~\cite{fan2008liblinear,pedregosa2011scikit}. Moreover, even when the user has a choice of objective function, the same function is used for every block~\cite{lightning_2016}.

In this work, we propose a generic coordinate descent method for the resolution of the convex optimization problem
\begin{equation}
\label{eq:pb_to_solve}
\min_{x \in \mathbb R^N} \revisionone{\frac 12 x^\top Q x} + \sum_{j=1}^J c^f_j f_j(A^f_j x - b^f_j) + \sum_{i=1}^I c^g_i g_i(D^g_i x^{(i)} - b^g_i) + \sum_{l=1}^L c^h_l h_l(A^h_l x - b^h_l) \;.
\end{equation}
We shall call $f_j$, $g_i$ and $h_l$ atom functions. Each of them may be different. We will assume that $f_j$'s are differentiable and convex, $g_i$'s and $h_l$'s are proximal-friendly convex functions. 
\revisionone{As before 
$A^f_j \in \mathbb R^{M^f_j \times N}$ and $A^h_l \in \mathbb R^{M^h_l \times N}$ are matrices, $D_i^g$ is a multiple of the identity matrix of size $N_i$, $b^f_j \in \mathbb R^{M^f_j}$, $b^g_i \in \mathbb R^{N_i}$ and $b^h_l \in \mathbb R^{M^h_l}$ are vectors, $c^f_j$, $c^g_i$ and $c^h_l$ are positive real numbers, $Q$ is a $N \times N$ positive semi-definite matrix.}

The algorithm we implemented is described in~\cite{fercoq2015coordinate}
and can be downloaded on \url{https://bitbucket.org/ofercoq/cd_solver}.
\revisionone{The present paper focuses on important implementation details about residual updates and dual variable duplication. The novelty of our code is that it} allows a generic treatment of these algorithmic steps and includes a modelling interface in Python for the definition of the optimization problem.
Note that unlike most coordinate descent implementations, it can deal with nonseparable nonsmooth objectives and linear constraints.

\revisionone{
\section{Literature review on coordinate descent methods}

A thorough review on coordinate descent is beyond the scope of this paper. 
We shall refer the interested reader to the review papers \cite{wright2015coordinate} and \cite{shi2016primer}.
Instead, for selected papers dealing with smooth functions, separable non-smooth functions or non-separable non-smooth function, we list their main features. We also quickly review what has been done for non-convex functions.
We sort papers in order of publication except when there are an explicit dependency between a paper and a follow-up.

\subsection{Smooth functions}

Smooth objectives are a natural starting point for algorithmic innovations. 
The optimization problem at stake writes
\[
\min_{x \in \mathbb R^N} f(x)
\]
where $f$ is a convex differentiable function with Lipschitz-continuous partial derivatives.

In Table~\ref{tab:smooth}, we compare papers that introduced important improvements to coordinate descent methods. 
We shall in particular stress the seminal paper by Tseng and Yun \cite{Tseng:CGDM:Nonsmooth}.
It features coordinate gradient steps instead of exact minimization. Indeed a 
coordinate gradient steps gives similar iteration complexity both in theory and in practice for a much cheaper iteration cost. Moreover, this opened the door for many innovations: blocks of coordinates and the use of proximal operators were developed in the same paper. Another crucial step was made in~\cite{nesterov2012efficiency}: Nesterov showed how randomization can help finding finite-time complexity bounds and proposed an accelerated version of coordinate descent. He also proposed to use a non-uniform sampling of coordinates depending on the coordinate-wise Lipschitz constants.

\begin{table}[htbp]
	\begin{center}
		{
			\begin{tabular}{| l || c | c | c | c | c | c || c |}
				\hline
				{\bf Paper} & \!\!{\bf Rate}\!\! & \!\!{\bf Rand}\!\! & \!\!{\bf Grad}\!\! & \!\!{\bf Blck}\!\! &  \!{\bf Par}\!  &  \!{\bf Acc }\! &
				{\bf  Notable feature} \\
				\hline
				Seidel '74 \cite{seidel1874ueber} & $\times$ & N & $\times$ & $\times$ & $\times$ & $\times$ & quadratic \\
				Warga '63 \cite{warga1963minimizing} & $\times$ & N & $\times$ & $\times$ & $\times$ & $\times$ & strictly  convex\\
				Luo \& Tseng '92 \cite{luo1992convergence} & \!\!asymp\!\! & N & $\times$ & $\times$ & $\times$ & $\times$ & rate for weakly convex \\
				Leventhal \& Lewis '08 \cite{Leventhal:2008:RMLC} 
				& \cmark & Y & $\times$ & $\times$ & $\times$ & $\times$ & quadratic $f$
				\\
				Tseng \& Yun '09 \cite{Tseng:CGDM:Nonsmooth} & \!\!asymp\!\! & N & \cmark & \cmark & $\times$ & $\times$ & \!\!line search, proximal operator\!\!\\
				Nesterov '12 \cite{nesterov2012efficiency} 
				& \cmark & Y & \cmark & $\times$ & $\times$ & \!\!not eff\!\! & 1st acc \& 1st non-uniform
				\\
				Beck \& Tetruashvili '13 \cite{beck2013convergence} &
				\cmark & N & \cmark & \cmark & $\times$ & \!\!not eff\!\! & 
				\!\!finite time analysis cyclic CD\!\!\\
				Lin \& Xiao '13 \cite{luxiao2013complexity}
				 & \cmark & Y & \cmark & \cmark & $\times$ & \!\!not eff\!\! & improvements on \cite{nesterov2012efficiency,RT:UCDC}
				\\
				Lee \& Sidford '13  \cite{lee2013efficient} 
			    & \cmark & Y & \cmark & $\times$ & $\times$ & \cmark &  1st efficient accelerated 
				\\
				Liu et al '13 \cite{Wright:2013-async_par_cdm}  
				& \cmark & Y & \cmark & $\times$ & \cmark & $\times$ & 1st asynchronous
				\\		
				Glasmachers \!\&\! Dogan \!'13\! \cite{glasmachers2013accelerated}\!\!\! & $\times$ & Y & \cmark & $\times$ & $\times$ & $\times$ & heuristic sampling \\
				Richt\'{a}rik \& Tak\'{a}\v{c} '16 \cite{RT:2013optimal} 
				& \cmark & Y & \cmark & $\times$  & \cmark  & $\times$ &  1st arbitrary sampling
				\\
				Allen-Zhu et al '16 \cite{allenzhu2016even} & \cmark & Y &\cmark & \cmark & $\times$ & \cmark & non-uniform sampling 				\\
				Sun et al '17 \cite{sun2017asynchronous} & \cmark & Y & \cmark & \cmark & \cmark & $\times$ & better asynchrony than \cite{Wright:2013-async_par_cdm} \\
				\hline
			\end{tabular}
		}
	\end{center}
	\caption{Selected papers for the minimization of smooth functions. {\bf Rate}: we check whether the paper proves convergence ($\times$), an asymptotic rate (asymp) or a finite time iteration complexity (\cmark). {\bf Rand}: deterministic (N) or randomized (Y) selection of coordinates. {\bf Grad}: exact minimization ($\times$) or use of partial derivatives (\cmark).
	{\bf Blck}: the paper considers 1D coordinates ($\times$) or blocks of coordinates (\cmark).
	{\bf Par}: Algorithm designed for parallel computing (\cmark).
	{\bf Acc}: no momentum ($\times$), accelerated but not efficient in practice (not eff), accelerated algorithm (\cmark)
	}
	\label{tab:smooth}
\end{table}

\subsection{Separable non-smooth functions}

A large literature has been devoted to composite optimization problems with separable non-smooth functions:
\[
\min_{x \in \mathbb R^N} f(x) + \sum_{i=1}^n g_i(x^{(i)})
\]
where $f$ is a convex differentiable function with Lipschitz-continuous partial derivatives and for all $i$, $g_i$ is a convex function whose proximal operator is easy to compute. Indeed, regularized expected risk minimization problems often fit into this framework and this made the success of coordinate descent methods for machine learning applications. 

Some papers study $\min_{x \in \mathbb R^p} f(x) + \sum_{i=1}^n g_i((A x)^{(i)})$, where $f$ is strongly convex and apply coordinate descent to a dual problem written as
\[
\min_{y \in \mathbb R^N} f^*(-A^\top y) + \sum_{i=1}^n g^*_i(y^{(i)}) \;.
\]
One of the challenges of these works is to show that even though we are solving the dual problem, one can still recover rates for a sequence minimizing the primal objective.

We present our selection of papers devoted to this type of problems in Table~\ref{tab:separable}.

\begin{table}[htbp]
	\begin{center}
		{
			\begin{tabular}{| l ||  c | c | c | c || c |}
				\hline
				{\bf Paper} & {\bf Prx} &  {\bf Par}  &  {\bf Acc } & {\bf Dual } &
				{\bf  Notable feature} \\
				\hline
				Tseng \& Yun '09 \cite{Tseng:CGDM:Nonsmooth} & \cmark & $\times$ & $\times$ & $\times$ & 1st prox, line search, deterministic\\		
				S-Shwartz \& Tewari '09 \cite{ShalevTewari09} 
				& $\ell_1$ & $\times$ & $\times$ & $\times$ & 1st $\ell_1$-regularized  w finite time bound
				\\
				Bradley et al '11 \cite{Bradley:PCD-paper} 
				& $\ell_1$ & \cmark & $\times$ & $\times$ & $\ell_1$-regularized parallel
								\\
				Richt\'{a}rik \& Tak\'{a}\v{c} '14 \cite{RT:UCDC} 
				& \cmark & $\times$ & $\times$ & $\times$ & 1st proximal with finite time bound
				\\				
				S-Shwartz \& Zhang '13 \cite{Proximal-dual-Coord-Ascent} 
				& \cmark & $\times$ & $\times$ & \cmark &  1st dual
				\\
				Richt\'{a}rik \& Tak\'{a}\v{c} '15 \cite{RT:PCDM}  
				& \cmark & \cmark& $\times$ & $\times$ & 1st general parallel
				\\
				Tak\'{a}\v{c} et al '13 \cite{minibatch-ICML2013} 
				& \cmark & \cmark & $\times$ & \cmark & 1st dual \& parallel
				\\
				S-Shwartz \& Zhang '14 \cite{shalev2013accelerated} 
				& \cmark & $\times$ & \cmark & \cmark & acceleration in the primal
				\\
				Yun '14 \cite{yun2014iteration} & \cmark & $\times$ & $\times$ & $\times$ &
				analysis of cyclic CD\\
				Fercoq \& Richt\'arik '15 \cite{FR:2013approx} & \cmark &  \cmark & \cmark & $\times$ & 1st proximal and accelerated
				\\
				Lin, Lu \& Xiao '14 \cite{lin2014accelerated} &\cmark & $\times$ & \cmark & $\times$ & prox \& accelerated on strong conv.\\
				Richt\'{a}rik \& Tak\'{a}\v{c} '16 \cite{RT:2013distributed} 
				& \cmark & \cmark & $\times$ & $\times$ & 1st distributed
				\\
				Fercoq et al '14 \cite{Hydra2} & \cmark &  \cmark & \cmark & $\times$ & distributed computation \\
				Lu \& Xiao '15 \cite{luxiao2013complexity} & \cmark & $\times$ & $\times$ & $\times$ & improved complexity over \cite{RT:UCDC,nesterov2012efficiency} 
				\\
				Li \& Lin '18 \cite{li2018complexity} & \cmark & $\times$ & \cmark & \cmark &  acceleration in the dual \\
				Fercoq \& Qu '18 \cite{fercoq2018restarting} & \cmark & $\times$ & \cmark & $\times$ & restart for obj with error bound \\
				\hline
			\end{tabular}
		}
	\end{center}
	\caption{An overview of selected papers proposing and analyzing the iteration complexity of coordinate descent methods for separable non-smooth objectives. {\bf Prx}: uses a proximal operator to deal with the non-smooth part of the objective. {\bf Par}: updates several blocks of coordinates in parallel.  {\bf Acc}: uses momentum to obtain an improved rate of convergence. {\bf Dual}: solves a dual problem but still proves rates in the primal (only relevant for weakly convex duals).}
	\label{tab:separable}
\end{table}

\subsection{Non-separable non-smooth functions}

Non-separable non-smooth objective functions are much more challenging to coordinate descent methods. 
One wishes to solve
\[
\min_{x \in \mathbb R^N} f(x) + g(x) + h(Ax)
\]
where $f$ is a convex differentiable function with Lipschitz-continuous partial derivatives, $g$ and $h$ are convex functions whose proximal operator are easy to compute and $A$ is a linear operator. 
Indeed, the linear operator introduces a coupling between the coordinates and a naive approach leads to a method that does not converge to a minimizer~\cite{Auslender1976}.
When $h = \iota_{\{b\}}$, the convex indicator function of the set $\{b\}$, we have equality constraints.

We present our selection of papers devoted to this type of problems in Table~\ref{tab:nonseparable}.

\begin{table}[htbp]
	\begin{center}
		{
			\begin{tabular}{| l || c | c | c  || c |}
				\hline
				{\bf Paper} & {\bf Rate} & {\bf Const} &  {\bf PD-CD}  &  
				{\bf  Notable feature} \\
				\hline
				Platt '99 \cite{platt199912} & $\times$ & \cmark & P  & for SVM 
				\\
				Tseng \& Yun '09 \cite{Tseng:CGDMLC:Nonsmooth} & \cmark & \cmark & P & adapts Gauss-Southwell rule \\
				Tao et al '12 \cite{tao2012stochastic} 
				& \cmark & $\times$ & P &  uses averages of subgradients
				\\
				Necoara et al '12 \cite{Necoara:Coupled} 
				& \cmark & \cmark & P & 2-coordinate descent 
				\\
				Nesterov '12 \cite{Nesterov-Subgrad-Huge} & \cmark & $\times$ & P & uses subgradients				\\
				Necoara \& Clipici '13 \cite{Necoara:parallelCDM-MPC} 
				& \cmark & \cmark & P& coupled constraints 
				\\
				Combettes \& Pesquet '14 \cite{combettes2014stochastic}
				& $\times$ & \cmark & \cmark & 1st PD-CD, short step sizes\\
				Bianchi et al '14 \cite{bianchi2014stochastic} & $\times$ & \cmark & \cmark& distributed optimization \\
				Hong et al '14 \cite{hong2014block} & $\times$ & \cmark & $\times$ & updates all dual variables\\
				Fercoq \& Richt\'{a}rik '17 \cite{FR:spcdm} 
				& \cmark & $\times$ & P & uses smoothing
				\\
				Alacaoglu et al '17 \cite{alacaoglu2017smooth}  & \cmark & \cmark & \cmark & 1st PD-CD w rate for constraints\\
				Xu \& Zhang '18 \cite{xu2018accelerated}& \cmark & \cmark  & $\times$& better rate than \cite{gao2016randomized}\\
				Chambolle et al '18 \cite{chambolle2017stochastic}  & \cmark & \cmark  & $\times$ & updates all primal variables\\
				Fercoq \& Bianchi '19 \cite{fercoq2015coordinate}& \cmark & \cmark  & \cmark & 1st PD-CD w long step sizes \\
				Gao et al '19 \cite{gao2016randomized} & \cmark & \cmark & $\times$ & 1st primal-dual w rate for constraints \\				
				Latafat et al '19 \cite{latafat2019new}  & \cmark & \cmark  & \cmark & linear conv w growth condition \\
	\hline
\end{tabular}
}
\end{center}
\caption{Selected papers for the minimization of non-smooth non-separable functions.
{\bf Rate}: does the paper prove rates? {\bf Const}: can the method solve problems with linear equality constraints? {\bf PD-CD}: the method is purely primal (P), the method updates some primal variables but all the Lagrange multipliers or some dual variables but all the primal variables ($\times$), the method updates some primal and some dual variables at each iteration (\cmark).
}
\label{tab:nonseparable}
\end{table}

\subsection{Non-convex functions}

The goal here is to find a local minimum to the problem
\[
\locmin_{x \in \mathbb R^N} f(x) + \sum_{i=1}^n g_i(x^{(i)})
\]
without any assumption on the convexity of $f$ nor $g_i$. The function $f$ should be continuously differentiable and the proximal operator of each function $g_i$ should be easily computable. Note that in the non-convex setting, the proximal operator may be set-valued.

\begin{table}[htbp]
	\begin{center}
		{
			\begin{tabular}{| l || c | c | c || c |}
				\hline
				{\bf Paper} & {\bf Conv} & {\bf Smth} & {\bf Nsmth} &
				{\bf  Notable feature} \\
	\hline
Grippo \& Sciandrone '00 \cite{grippo2000convergence}\!\!& $\times$ & \cmark & $\times$ & 2-block or coordinate-wise quasiconvex \\
Tseng \& Yun '09 \cite{Tseng:CGDM:Nonsmooth} & \cmark & \cmark & $\times$ & convergence under error bound  \\
Hsieh \& Dhillon '11 \cite{hsieh2011fast} & $\times$ & \cmark & $\times$ & non-negative matrix factorization\\
Breheny \& Huang '11 \cite{breheny2011coordinate} & $\times$ & $\times$ & \cmark & regularized least squares\\
Mazumder et al '11 \cite{mazumder2011sparsenet} & $\times$ & $\times$ & \cmark & regularized least squares \\
Razaviyayn et al '13 \cite{razaviyayn2013unified} & $\times$ & \cmark & \cmark & requires uniqueness of the prox\\
Xu \& Yin '13 \cite{xu2013block} & \cmark & \cmark & $\times$ & multiconvex \\
Lu \& Xiao '13 \cite{lu2013randomized} & $\times$ & \cmark & \cmark & random sampling \\
Patrascu \& Necoara '15 \cite{patrascu2015efficient} & \cmark & \cmark & $\times$ & randomized + 1 linear constraint\\ 
Xu \& Yin '17 \cite{xu2017globally} & \cmark & \cmark & \cmark & convergence under KL \\
	\hline 
\end{tabular}
}
\end{center}
\caption{Selected papers for the local minimization of nonconvex problems. {\bf Conv}: study only limit points ($\times$) or convergence of the sequence proved (\cmark). {\bf Smth}: can deal with non-convex smooth functions. {\bf Nsmth}: can deal with non-convex and non-smooth functions.}
\end{table}
}

\section{Description of the Algorithm}

\subsection{General scheme}

The algorithm we implemented is a coordinate descent primal-dual method developed in~\cite{fercoq2015coordinate}. Let introduce the notation $F(x) = \frac 12 x^\top Q x +  \sum_{j=1}^J c^f_j f_j(A^f_j x - b^f_j)$, $G(x) = \sum_{i=1}^I c^g_i g_i(D^g_i x^{(i)} - b^g_i) $, $H(z) = \sum_{l=1}^L c^h_l h_l(z^{(l)} - b^h_l)$, $\mathcal J(i) = \{j \;:\; A^h_{j,i} \neq 0\}$, $\mathcal I(j) = \{i \;:\; A^h_{j,i} \neq 0\}$, $m_j = |\mathcal I(j)|$ and $\rho(A)$ the spectral radius of matrix $A$. 
\revisionone{We shall also denote $\mathcal J^f(i) = \{j \;:\; A^f_{j,i} \neq 0\}$, $\mathcal J^Q(i) = \{j \;:\; Q_{j,i} \neq 0\}$, $A^f \in \mathbb R^{\sum_j M^f_j \times N}$ the matrix which stacks the matrices $(A^f_j)_{1 \leq j \leq J}$ and $A^h$ the matrix which stacks the matrices $(A^h_l)_{1 \leq l \leq L}$.}
The algorithm writes then as Algorithm~\ref{algo:main}.

\begin{algorithm}[ht]
\caption{Coordinate-descent primal-dual algorithm with duplicated variables (PD-CD)}
\label{algo:main}
  \noindent {\bf Input}: Differentiable function $F: \mathbb R^N \to \mathbb R$, matrix $A^h \in \mathbb R^{M^h \times N}$, functions $G$ and $H$ whose proximal operators are available.
  
  \noindent {\bf Initialization}: Choose $x_0\in \mathbb R^N$, $\bs y_0\in \mathbb R^{\mathrm{nnz}(A^h)}$.
\revisionone{Denote  $\mathcal J(i) = \{j \;:\; A^h_{j,i} \neq 0\}$, $\mathcal I(j) = \{i \;:\; A^h_{j,i} \neq 0\}$, $m_j = |\mathcal I(j)|$ and $\rho(A)$ the spectral radius of matrix $A$. }
Choose step sizes $\tau \in \mathbb R_+^I$ and $\sigma \in \mathbb R_+^L$
such that $\forall i \in \{1, \ldots I\}$,
\begin{equation}
\label{eq:tau_algo2}
\tau_i < \frac 1{\beta_i+\rho\left(\sum_{j\in \mathcal J(i)}m_j\sigma_j (A^h)_{j,i}^\top A^h_{j,i}\right)}\,.
\end{equation}
 \\ For all $i \in \{1,\dots,I\}$, set $w_0^{(i)} = \sum_{j\in
\mathcal J(i)}(A^h)_{j,i}^\top\,\bs y_0^{(j)}(i)$.\\ For all $j \in \{1,\dots,J\}$, set
$z_0^{(j)} = \frac 1{m_j}\sum_{i\in \mathcal I(j)}\bs y_0^{(j)}(i)$.

  \noindent {\bf Iteration $k$}: Define:
  \begin{align*}
    \overline y_{k+1} &= \prox_{\sigma,H^\star}\big(z_k+D(\sigma) A^h x_k\big) \\
    \overline x_{k+1} &= \prox_{\tau,G}\Big(x_k-D(\tau)\left(\nabla F(x_k)+2(A^h)^\top \overline {y}_{k+1}- w_k\right)\Big)\,.
  \end{align*}
  \hspace{2em} For $i=i_{k+1} \sim U(\{1, \ldots, I \})$ and for each $j \in \mathcal J(i_{k+1})$, update:
  \begin{align*}
    &{x}_{k+1}^{(i)} = \overline {x}_{k+1}^{(i)} \\
    &{\bs y}_{k+1}^{(j)}(i) = \overline {y}_{k+1}^{(j)}\\
    &w_{k+1}^{(i)} = w_{k}^{(i)} + \sum_{j \in J(i)}(A^h)_{j,i}^\top\,({\bs y}_{k+1}^{(j)}(i)-\bs y_k^{(j)}(i))\\
    &z_{k+1}^{(j)} = z_{k}^{(j)} + \frac 1{m_j} ({\bs y}_{k+1}^{(j)}(i)
    - \bs y_k^{(j)}(i)) \,.
  \end{align*}
  \hspace{2em}   Otherwise, set ${x}_{k+1}^{(i')}=x_k^{(i')}$,
  $w_{k+1}^{(i')}=w_k^{(i')}$, $z_{k+1}^{(j')}=z_k^{(j')}$ and ${\bs y}_{k+1}^{(j')}(i')={\bs y}_{k}^{(j')}(i')$.
\end{algorithm}

We will denote $U_1, \ldots, U_I$ the columns of the identity matrix
corresponding to the blocks of $x = (x^{(1)}, \ldots, x^{(I)})$, so that $U_i x^{(i)} \in \mathbb R^N$  and
$V_1, \ldots, V_J$ the columns of the identity matrix
corresponding to the blocks of $A^f x - b^f = (A^f_1 x - b^f_1, \ldots, A^f_J x - b^f_J)$.

\subsection{Computation of partial derivatives}

For simplicity of implementation, we are assuming that $G$ is separable 
and the blocks of variable will follow the block structure of $G$.
This implies in particular that at each iteration, only $\nabla_i F(x_k)$ 
needs to be computed. 
This partial derivative needs to be calculated efficiently because it needs to be performed at each iteration of the algorithm. 
We now describe the efficient residual update method, which is classically used in coordinate descent implementations~\cite{nesterov2012efficiency}.

Denote $r^{f,x}_k = A^f x_k - b^f$ and $r^{Q,x}_k = Q x_k$. By the chain rule, we have 
\[
\nabla_i F(x_k) = \sum_{j=1}^J c^f_j (A^f)^\top_{j,i} \nabla f_j(A^f_j x_k - b^f_j) + Q x_k= \sum_{j \in \mathcal J^f(i)} c^f_j (A^f)^\top_{j,i} \nabla f_j((r^{f,x}_k)_j) + \sum_{j \in \mathcal J^Q(i)} U_i r^{Q,x}_k
\]
If $r^{f,x}_k$ and $r^{Q,x}_k$ are pre-computed, only $O(|\mathcal J^f(i)| + |\mathcal J^Q(i)|)$ operations are needed.

For an efficient implementation, we will update the residuals $r^{f,x}_k$ as follows, using the fact that only the coordinate block $i_{k+1}$ is updated:
\begin{multline*}
r^{f,x}_{k+1} = A^f x_{k+1} - b^f = A^f \big(x_{k} + U_{i_{k+1}}(x_{k+1}^{(i_{k+1})} - x_k^{(i_{k+1})})\big) - b^f = r^{f,x}_k + A^f U_{i_{k+1}}(x_{k+1}^{(i_{k+1})} - x_k^{(i_{k+1})})\\
= r^{f,x}_k + \sum_{j \in \mathcal J^f(i_{k+1})} V_j A^f_{j,i_{k+1}} (x_{k+1}^{(i_{k+1})} - x_k^{(i_{k+1})}) 
\end{multline*}
Hence, updating $r^{f,x}_{k+1}$ also requires only $O(|\mathcal J^f(i_{k+1})|)$ iterations.

Similarly, updating the residuals $r^{Q,x}_k$, $r^{h,x}_k = A^h x_k - b^h$,
$w_k$ and 
$z_k$ 
can be done in $O(|\mathcal J^Q(i_{k+1})|)$ and $O(|\mathcal J(i_{k+1})|)$ operations.

\revisionone{Although this technique is well known, it is not trivial how to write it in a generic fashion, since residual updates are needed at each iteration and should be written in a compiled language. We coded the residual update using abstract atom functions in order to achieve this goal.}

\subsection{Computation of proximal operators using atom functions}

Another major step in the method is the computation of the $i^{\text{th}}$ coordinate of $\prox_{\tau, G}(x)$ for a given $x \in \mathbb R^N$.

As $D^g$ is assumed to be diagonal, $G$ is separable. Hence, by the change of variable $\bar z = D^g_i \bar x - b^g_i$,
\begin{align*}
(\prox_{\tau, G}(x))_i & = \arg \min_{\bar x \in \mathbb R^{N_i}} c_i^g g_i(D^g_i \bar x - b^g_i) + \frac{1}{2 \tau_i} \norm{\bar x - x^{(i)}}^2\\
& = (D_i^g)^{-1} \Big(b^g_i+\arg \min_{\bar z \in \mathbb R^{N_i}} c_i^g g_i(\bar z) + \frac{1}{2 \tau_i} \norm{(D_i^g)^{-1} (b^g_i+\bar z) - x^{(i)}}^2\Big)\\
& = (D_i^g)^{-1} \Big(b^g_i+\arg \min_{\bar z \in \mathbb R^{N_i}} g_i(\bar z) + \frac{1}{2 c_i^g (D_i^g)^2 \tau_i } \norm{\bar z - (D^g_i x^{(i)} - b^g_i)}^2\Big)\\
& = (D_i^g)^{-1} \Big(b^g_i + \prox_{c_i^g (D_i^g)^2 \tau_i g_i}(D_i^g x^{(i)} - b_i^g)\Big)
\end{align*}
where we used the abuse of notation that $D_i^g$ is either the scaled identity matrix or any of its diagonal elements.
This derivation shows that to compute $(\prox_{\tau, G}(x))_i$ we only need
linear algebra and the proximal operator of the atom function $g_i$.

We can similarly compute $\prox H$. To compute $\prox_{\sigma, H^\star}$,
we use Moreau's formula:
\[
\prox_{\sigma, H^\star}(z) = z - D(\sigma) \prox_{\sigma^{-1}, H}(D(\sigma)^{-1}z) 
\]

\subsection{Duplication of dual variables}

Algorithm~\ref{algo:main} maintains duplicated dual variables $\bs y_k \in \mathbb R^{\mathrm{nnz}(A^h)}$
as well as averaged dual variables $z_k \in \mathbb R^{M^h}$ where $M^h = \sum_{l=1}^L M^h_l$ and $A^h_{l,i}$ is of size $M_l^h \times N_i$.
The sets $\mathcal J(i)$ for all $i$ are given by the sparse column format representation of $A^h$. Yet, for all $i$, we need to construct the set of
indices of $\bs y_{k+1}$ that need to be updated. This is the table
\texttt{dual\_vars\_to\_update} in the code.
Moreover, as $H$ is not separable in general, in order to compute $\bar y_{k+1}^j$, for $j \in \mathcal J(i_{k+1})$, we need to determine the set of dual indices $j'$ that belong to 
the same block as $j$ with respect to the block decomposition of $H$.
This is the purpose of the tables \texttt{inv\_blocks\_h} and \texttt{blocks\_h}. 

\revisionone{The procedure allows us to only compute the entries of $\bar y_{k+1}$ that are required for the update of $\boldsymbol y_k$.}

\section{Code structure}

\begin{figure}[htb]
\centering
\includegraphics[width=0.7\linewidth]{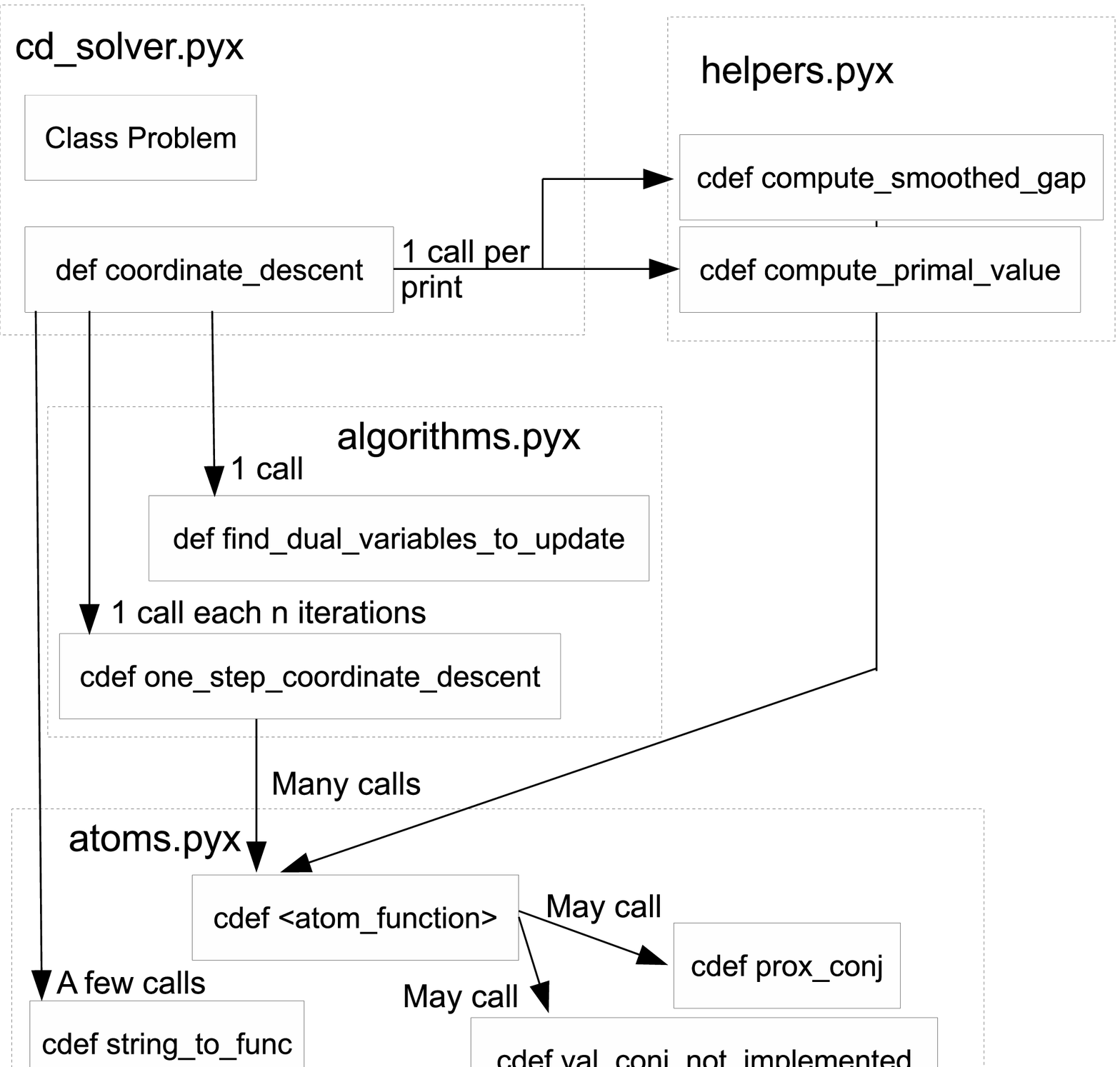}
\caption{Code structure when using no screening and no inertial acceleration}
\label{fig:code_structure}
\end{figure}

The code is organized in nine files. The main file is \texttt{cd\_solver.pyx}. It contains the Python callable and the data structure for the problem definition. 
The other files are \texttt{atoms.pyx/pxd}, \texttt{algorithm.pyx/pxd}, \texttt{helpers.pyx/pxd} and \texttt{screening.pyx/pxd}. They contain the definition of the atom functions, the algorithms
and the functions for computing the objective value. In Figure~\ref{fig:code_structure}, we show for each subfunction, in which function it is used. The user needs to call the Python class \texttt{Problem}
and the Python function \texttt{coordinate\_descent}.
Atom functions can be added by the user without modifying the main algorithm.

All tables are defined using Numpy's array constructor in the \texttt{coordinate\_descent} function. The main loop of coordinate descent and the atom functions are pre-compiled for efficiency.

\section{Atom functions}

The code allows us to define atom functions independently of the 
coordinate descent algorithm. 
As an example, we provide in Figure~\ref{fig:square_atom} the code for the square function atom.

\begin{figure}[htbp]
\begin{verbatim}
cdef DOUBLE square(DOUBLE[:] x, DOUBLE[:] buff, int nb_coord, MODE mode, 
        DOUBLE prox_param, DOUBLE prox_param2) nogil:
    # Function x -> x**2
    cdef int i
    cdef DOUBLE val = 0.
    if mode == GRAD:
        for i in range(nb_coord):
            buff[i] = 2. * x[i]
        return buff[0]
    elif mode == PROX:
        for i in range(nb_coord):
            buff[i] = x[i] / (1. + 2. * prox_param)
        return buff[0]
    elif mode == PROX_CONJ:
        return prox_conj(square, x, buff, nb_coord, prox_param, prox_param2)
    elif mode == LIPSCHITZ:
        buff[0] = 2.
        return buff[0]
    elif mode == VAL_CONJ:
        return val_conj_not_implemented(square, x, buff, nb_coord)
    else:  # mode == VAL
        for i in range(nb_coord):
            val += x[i] * x[i]
        return val
\end{verbatim}
\vspace{-2ex}
\caption{Code for the square function atom}
\label{fig:square_atom}
\end{figure}

As inputs, it gets \texttt{x} (an array of numbers which is the point where the operation takes place), \texttt{buff} (the buffer for vectorial outputs), \texttt{nb\_coord} (is the size of \texttt{x}), \texttt{mode}, \texttt{prox\_param} and \texttt{prox\_param2} (numbers which are needed when computing the proximal operator).
The input \texttt{mode} can be:
\begin{itemize}
\item \texttt{GRAD} in order to compute the gradient. 
\item \texttt{PROX} to compute the proximal operator.
\item \texttt{PROX\_CONJ} uses Moreau's formula to compute the proximal operator of the conjugate function.
\item \texttt{LIPSCHITZ} to return the Lipschitz constant of the gradient. 
\item \texttt{VAL\_CONJ} to return the value of the 
conjugate function. As this mode is used only by \texttt{compute\_smoothed\_gap} for printing purposes, its implementation is optional and can be approximated using the helper function \texttt{val\_conj\_not\_implemented}. \revisionone{Indeed, for a small $\epsilon > 0$,
$h^*(y) = \sup_{z} \langle z, y \rangle - h(z) \approx  \sup_{z} \langle z, y \rangle - h(z) - \frac{\epsilon}{2} \norm{z}^2 = \langle p, y \rangle - h(p) - \frac{\epsilon}{2} \norm{p}^2$,
where $p = \prox_{h / \epsilon}(y/\epsilon)$.}
\item \texttt{VAL} to return the value of the function.
\end{itemize} 

Some functions naturally require multi-dimensional inputs, like $\norm{\cdot}_2$ or the log-sum-exp function. For consistency, we define all the atoms with multi-dimensional inputs: for an atom function $f_0:\mathbb R \to \mathbb R$, we extend it to an
atom function $f:\mathbb R^{N_i} \to \mathbb R$ by $f(x) = \sum_{l=1}^{N_i}
f_0(x_l)$. 

For efficiency purposes, we are bypassing the square atom function when computing a gradient and implemented it directly in the algorithm.

\section{Modelling language}

In order to use the code in all its generality, we defined a modelling
language that can be used to define the optimization problem we want to solve~\eqref{eq:pb_to_solve}. 

The user defines a problem using the class \texttt{Problem}.
Its arguments can be:
\begin{itemize}
\item \texttt{N} the number of variables, \texttt{blocks} the blocks of coordinates coded in the same fashion as the \texttt{indptr} index of sparse matrices (default [0, 1, \ldots, \texttt{N}]), \texttt{x\_init} the initial primal point (default 0) and \texttt{y\_init} the initial duplicated dual variable (default 0)

\item Lists of strings \texttt{f}, \texttt{g} and \texttt{h} that code for the atom functions used. The function \texttt{string\_to\_func} is responsible for 
linking the atom function that corresponds to the string. Our convention
is that the string code is exactly the name of the function in \texttt{atoms.pyx}. The size of the input of each atom function is
defined in \texttt{blocks\_f}, \texttt{blocks} and \texttt{blocks\_h}.
The function strings \texttt{f}, \texttt{g} or \texttt{h} may be absent,
which means that the function does not appear in the problem to solve.

\item Arrays and matrices \texttt{cf}, \texttt{Af}, \texttt{bf}, \texttt{cg}, \texttt{Dg}, \texttt{bg}, \texttt{ch}, \texttt{Ah}, \texttt{bh}, \texttt{Q}.
The class initiator transforms matrices into the sparse column format and checks whether \texttt{Dg} is diagonal.
\end{itemize}
For instance, in order to solve the Lasso problem, $\min_x \frac 12 \norm{Ax - b}^2_2 + \lambda \norm{x}_1$, one can 
type

\noindent \texttt{pb\_lasso = cd\_solver.Problem(N=A.shape[1],
        f=["square"]*A.shape[0],
        Af=A,}
        
\texttt{bf=b,
        cf=[0.5]*A.shape[0],
        g=["abs"]*A.shape[1],
        cg=[lambda]*A.shape[1])}

\noindent \texttt{cd\_solver.coordinate\_descent(pb\_lasso)}

\section{Extensions}

\subsection{Non-uniform probabilities}

We added the following feature for an improved efficiency. Under the argument \texttt{sampling='kink\_half'}, the algorithms periodically detects the set of blocks $I_{\mathrm{kink}}$ such that $i \in I_{\mathrm{kink}}$ if $x^{(i)}$ is at a kink of $g_i$. 
Then, block $i$ is selected with probability law
\[
\mathbb P(i_{k+1} = i) = \begin{cases}
\frac{1}{n} & \text{ if } |I_{\mathrm{kink}}| = n \\
\frac{1}{2n} & \text{ if } |I_{\mathrm{kink}}| < n \text{ and } i \in I_{\mathrm{kink}} \\ \frac{1}{2n} + \frac{1}{2 (n - |I_{\mathrm{kink}}|)}  & \text{ if } |I_{\mathrm{kink}}| < n \text{ and } i \not \in I_{\mathrm{kink}}
\end{cases}
\]
The rationale for this probability law is that blocks at kinks are likely to incur no move when we try to update them. We thus put more computational power for non-kinks. On the other hand, we still keep an update probability weight of at least $\frac{1}{2n}$ for each block, so even in unfavourable cases, we should not observe too much degradation in the performance as compared to the uniform law.

\subsection{Acceleration}

We also coded accelerated coordinate descent~\cite{FR:2013approx}, as well as its restarted~\cite{fercoq2018restarting}
and primal-dual~\cite{alacaoglu2017smooth} variants.
The algorithm is given in Algorithm~\ref{algo:accel}. As before, $\bar y_{k+1}$ and $\bar x_{k+1}$ should not be computed: only the relevant coordinates should be computed. 

\begin{algorithm}[htbp]
\caption{Smooth, accelerate, randomize the Coordinate Descent (APPROX/SMART-CD)}
\label{algo:accel}
  \noindent {\bf Input}: Differentiable function $F: \mathbb R^N \to \mathbb R$, matrix $A^h \in \mathbb R^{M^h \times N}$, functions $G$ and $H$ whose proximal operators are available.
  
  \noindent {\bf Initialization}: Choose $x_0\in \mathbb R^N$, $\dot y_0\in \mathbb R^{M^h}$.
Choose $\gamma_1>0$ and denote $B_0^i = \beta_i + \frac{\rho((A^h_{:,i})^\top A^h_{:,i})}{\gamma_1}$.

Set $s=0$, $\theta_0 = \frac{1}{n}$, $c_0 =1$, $\hat x_0 = 0 \in \mathbb R^N$ and $\tilde x_0 = x_0$.

  \noindent {\bf Iteration $k$}: Define:
  \begin{align*}
    \overline y_{k+1} &= \prox_{\gamma^{-1}_{k+1},H^\star}\big(\dot y_s + D(\gamma_{k+1})^{-1}( c_k A^h \hat x_k + A^h \tilde x_k )\big) \\
    \overline x_{k+1} &= \prox_{\frac{\theta_0}{\theta_k} B_k^{-1},G}\Big(\tilde x_k-\frac{\theta_0}{\theta_k}D(B_k)^{-1}\left(\nabla F(c_k \hat x_k + \tilde x_k)+(A^h)^\top \overline{y}_{k+1}\right)\Big)\,. 
  \end{align*}
  \hspace{2em} For $i=i_{k+1} \sim U(\{1, \ldots, I \})$, update:
  \begin{align*}
    &\tilde {x}_{k+1}^{(i)} = \overline {x}_{k+1}^{(i)}\\
    & \hat x_{k+1}^{(i)} = \hat x_{k}^{(i)} - \frac{1-\theta_k/\theta_0}{c_k} (\tilde {x}_{k+1}^{(i)} - \tilde {x}_{k}^{(i)})
  \end{align*}
  \hspace{2em}   Otherwise, set ${x}_{k+1}^{(i')}=x_k^{(i')}$.
  
  \hspace{2em} Compute $\theta_{k+1} \in (0,1)$ as the unique positive root of
  \[
  \begin{cases}
  \theta^3 + \theta^2 + \theta_k^2 \theta - \theta_k^2 = 0 & \text{ if } h \neq 0 \\
  \phantom{\theta^3~~~~}\theta^2 + \theta_k^2 \theta - \theta_k^2 = 0 & \text{ if } h = 0 \\  
  \end{cases}
  \]
  
  \hspace{2em} Update $\gamma_{k+2} = \frac{\gamma_{k+1}}{1+\theta_{k+1}}$, $c_{k+1} = (1-\theta_{k+1}) c_k$ and $B_{k+1}^i = \beta_i + \frac{\rho((A^h_{:,i})^\top A^h_{:,i})}{\gamma_{k+2}}$ for all $i$.

  \noindent {\bf If } Restart($k$) is true:
  
 \hspace{2em} Set $x_{k+1} = \tilde x_{k+1} + c_{k} \hat x_{k+1}$.
 
 \hspace{2em} Set $\dot y_{s+1} = \overline y_{k+1}$ and $s \leftarrow s+1$.
 
 \hspace{2em} Reset $\hat x_{k+1} \leftarrow 0$, $\tilde x_{k+1} \leftarrow x_{k+1}$, $c_{k+1} \leftarrow 1$, $\theta_{k+1} \leftarrow \theta_0$, $\beta_{k+1} \leftarrow \beta_1$.
\end{algorithm}

The accelerated algorithms improve the worst case guarantee as explained in Table~\ref{tab:accel}:
\begin{table}[htbp]
\centering
\begin{tabular}{lcc}
					  &  $h = 0$   &  $h \neq 0$ \\
\hline
PD-CD Alg. \ref{algo:main}  & $O(1/k)$   & $O(1/\sqrt{k})$ \\
APPROX / SMART-CD Alg. \ref{algo:accel}    & $O(1/k^2)$ & $O(1/k)$\\
\hline
\end{tabular}
\caption{Convergence speed of the algorithms implemented}
\label{tab:accel}
\end{table}

However, accelerated algorithms do not take profit of regularity properties of the objective like strong convexity. Hence, they are not guaranteed to
be faster, even though restart may help.

\subsection{Variable screening}

The code includes the Gap Safe screening method presented in~\cite{eugene-phd}. Note that the method has been studied only for the case where $h=0$. Given a non-differentiability point $\hat x^{(i)}$ of the function $g_i$
where the subdifferential $\partial g_i(\hat x^{(i)})$ has a non-empty interior, 
a test is derived to check whether $\hat x^{(i)}$ is the $i^{\text{th}}$ variable of an optimal solution. If this is the case, one can set $x^{(i)} = \hat x^{(i)}$ and stop updating this variable. This may lead to a huge speed up in some cases. As the test relies on the computation of the duality gap, which has a nonnegligible cost, it is only performed from time to time.

\revisionone{In order to state Gap Safe screening in a general setting, we need the concept of polar of a support function.
Let $C$ be a convex set. The support function of $C$ is $\sigma_C$ defined by 
\[
\sigma_C(x) = \sup_{y \in C} \langle y, x \rangle \; .
\]
The polar to $\sigma_C$ is
\[
\sigma_C^\circ (\bar x) = \sup_{x:\sigma_C(x) \leq 1} \langle \bar x, x \rangle \;.
\]
In particular, if $\bar x \in C$, then $\sigma_C^\circ(\bar x) \leq 1$. Denote $f(z) = \sum_{j=1}^J c^f_j f_j(z^{(j)} - b^f_j)$, so that $(\nabla f(A^f x))_j = c^f_j \nabla f_j(A^f_j x - b_j^f)$, $G_i(x) = c^g_i g_i(D_i^g x^{(i)} - b_i^g)$, $x_\star$ a solution to the optimization problem~\eqref{eq:pb_to_solve} and suppose we have a set $\mathcal R$ such that $(\nabla f(A^f x_\star - b^f), Q x_\star) \in \mathcal R$. Gap Safe screening states that
\[
\max_{(\zeta, \omega) \in \mathcal R} \sigma^\circ_{\partial G_i(\hat x^{(i)})} ((A^f)_i^\top \zeta + \omega^{(i)}) < 1  \quad \Rightarrow \quad   x^{(i)}_\star = \hat x^{(i)} \;.
\]

Denote $z = \nabla f(A^f x)$ and $w = Qx$\footnote{We reuse the notation $w$ and $z$ here for the purpose of this section.}. We choose $\mathcal R$ as a sphere centered at 
\[
(\bar \zeta, \bar \omega) = \frac{(z, w)}{\max(1, \max_{1 \leq i \leq I}\big(\sigma_{\dom G_i^*}^\circ((A_i^f)^\top z+w^{(i)})\big))}
\]
and with radius
\[
r = \sqrt{\frac{2 \mathrm{Gap}(x, \bar \zeta, \bar \omega)}{L_{f,Q}}}
\]
where 
\[
L_{f, Q} = \max(\max_{1 \leq j \leq J} L(\nabla f_j), \max_{1 \leq i \leq I} \rho(Q_{i,i}))
\]
and
\[
\mathrm{Gap}(x, \bar \zeta, \bar \omega) = \frac 12 x^\top Q x + f(A^f x) + G(x) + G^*(-(A^f)^\top \bar \zeta - \bar \omega) + \frac 12 \bar \omega^\top Q^{\dagger} \bar \omega + f^*(\bar \zeta) \;.
\]
Note that as $\bar \omega$ is a rescaled version of $w = Qx$, we do not need to know $Q^{\dagger}$ in order to compute $\bar \omega^\top Q^{\dagger} \bar \omega$. It is proved in~\cite{eugene-phd} that this set $\mathcal R$ contains $(\nabla f(A^f x_\star - b^f), Q x_\star)$ for any optimal solutions $x_\star$ to the primal problem.
In the case where $G_i$ is a norm and $x^{(i)} = 0$, these expressions simplify since the $\sigma^\circ_{\partial G_i(0)} = \sigma^\circ_{\dom G_i^*}$ is nothing else than the dual norm associated to $G_i$.

For the estimation of $\max_{(\zeta, \omega) \in \mathcal R} \sigma^\circ_{\partial G_i(\hat x^{(i)})} ((A^f)_i^\top \zeta + \omega^{(i)})$, we use the fact that the polar of a support function is sublinear and positively homogeneous. Indeed, we have
\begin{align*}
\sigma^\circ_{\partial G_i(\hat x^{(i)})} &((A^f)_i^\top \zeta + \omega^{(i)}) = 
\sup_{x\;:\;\sigma_{\partial G_i(\hat x^{(i)})}(x) \leq 1} \langle (A^f)_i^\top \zeta + \omega^{(i)}, x \rangle =  \sup_{x\;:\;\sigma_{c^g_i D^g_i \partial g_i(D^g_i \hat x^{(i)} - b^g_i) \leq 1}} \langle (A^f)_i^\top \zeta + \omega^{(i)}, x \rangle \\
& = \frac{1}{c^g_i D^g_i}\sigma^\circ_{\partial g_i(D^g_i \hat x^{(i)} - b^g_i)} ((A^f)_i^\top \zeta+\omega^{(i)}) = \frac{1}{c^g_i D^g_i}\sigma^\circ_{\partial g_i(\hat x_{g_i})} ((A^f)_i^\top \zeta+\omega^{(i)}) \\
&\leq \frac{1}{c^g_i D^g_i}\sigma^\circ_{\partial g_i(\hat x_{g_i})} ((A^f)_i^\top \bar \zeta+ \bar\omega^{(i)}) + r \sup_{(u,v): \norm{(u,v)} = 1} \frac{1}{c^g_i D^g_i}\sigma^\circ_{\partial g_i(\hat x_{g_i})} ((A^f)_i^\top u + v^{(i)}) \;.
\end{align*}
Here $\hat x_{g_i} = D^g_i \hat x^{(i)} - b^g_i$ is a point where $\partial g_i(\hat x_{g_i})$ has a nonempty interior.
Some care should be taken when $\partial G_i(\hat x^{(i)})$ is unbounded, so that we first check whether $(A^f)_i^\top \zeta + \omega^{(i)} \in \dom \sigma^\circ_{\partial g_i(\hat x_{g_i})}$ for all $(\zeta, \omega) \in \mathcal R$.

Here also, the novelty lies in the genericity of the implementation.}

\section{Numerical validation}

\subsection{Performance}

In order to evaluate the performance of the implementation, we compare our implementation with a pure Python coordinate descent solver and code written for specific problems: Scikit learn's Lasso solver and Liblinear's SVM solver. We run the code on an
Intel Xeon CPU at 3.07GHz.

\begin{table}[htbp]
\centering
\begin{tabular}{lr}
Lasso & \\
\hline
Pure Python  &  308.76s \\
\texttt{cd\_solver}  & 0.43s \\
Scikit learn Lasso & 0.11s \\
\hline
\end{tabular}
\qquad \qquad 
\begin{tabular}{lr}
SVM & \\
\hline
Pure Python & 126.24s \\
\texttt{cd\_solver}  & 0.31s \\
Liblinear SVM  & 0.13s \\
\hline
\end{tabular}
\caption{Comparison of our code with a pure Python code and reference implementations for performing 100$N$ coordinate descent iterations for the Lasso problem on the Leukemia dataset with regularization parameter 
$\lambda = 0.1 \norm{(A^f)^\top b^f}_\infty$, and for 10$N$ coordinate descent iterations for the dual SVM problem on the RCV1 dataset with penalty parameter $C = 10$.}
\label{tab:time_comparison}
\end{table}

We can see on Table~\ref{tab:time_comparison} that our code is hundreds of times faster than the pure Python code. This is due to the compiled nature of our code, that does not suffer from the huge number of iterations required by
coordinate descent.
On the other hand, our code is about 4 times slower than state-of-the-art coordinate descent implementations designed for
a specific problem. We can see it in both examples we chose. This overhead is the price of genericity. 

We believe that, except for critical applications like Lasso or SVM, a 4 times speed-up does not justify writing
a new code from scratch, since a separate piece of code for each problem makes it difficult to maintain and to improve with future algorithmic advances.

\subsection{Genericity}

We tested our algorithm on the following problems:
\begin{itemize}
\item Lasso problem
\[
\min_{x \in \mathbb R^n} \frac 12 \norm{Ax - b}_2^2 + \lambda \norm{x}_1
\]
\item Binomial logistic regression
\[
\min_{x \in \mathbb R^n} \sum_{i=1}^m \log(1+\exp(b_i (Ax)_i)) + \frac{\lambda}{2} \norm{x}_2^2
\]
where $b_{i} \in \{-1, 1\}$ for all $i$.
\item Sparse binomial logistic regression
\[
\min_{x \in \mathbb R^n} \sum_{i=1}^m \log(1+\exp(b_i (Ax)_i)) + \lambda \norm{x}_1
\]

\item Dual SVM without intercept
\[
\min_{x \in \mathbb R^n} \frac{1}{2\alpha} \norm{A^\top D(b) x}^2_2 - e^\top x + \iota_{[0,1]^n}(x)
\]
where $\iota_{[0,1]^n}$ is the convex indicator function of the set $[0,1]^n$ and encodes the constraint $x \in [0,1]^n$.
\item Dual SVM with intercept
\[
\min_{x \in \mathbb R^n} \frac{1}{2\alpha} \norm{A^\top D(b) x}^2_2 - e^\top x + \iota_{[0,1]}(x) + \iota_{\{0\}}(b^\top x)
\]
\item Linearly constrained quadratic program
\[
\min_{x \in \mathbb R^n} \frac{1}{2} \norm{(A^f)^\top x - b^f}^2_2 + \iota_{\{0\}}(A^h x - b^h)
\]
\item Linear program
\[
\min_{x \in \mathbb R^n} c^\top x + \iota_{\mathbb R^n_+}(x) + \iota_{\mathbb R^m_-}(Ax - b)
\]
\item TV+$\ell_1$-regularized regression
\[
\min_{x \in \mathbb R^{n_1 n_2 n_3}} \frac 12 \norm{Ax - b}_2^2 + \alpha_1 \norm{Dx}_{2,1}
 + \alpha_2 \norm{x}_1\]
where $D$ is the discrete gradient operator and $\norm{y}_{2,1} = \sum_{i,j,k} \sqrt{\sum_{l=1}^3 y_{i,j,k,l}^2}$.
\item Sparse multinomial logistic regression
\[
\min_{x \in \mathbb R^{n \times q}} \sum_{i=1}^m \log\Big(\sum_{j=1}^q \exp\big(\sum_{l=1}^n A_{i,l} x_{l,j}\big)\Big) + \sum_{i=1}^n \sum_{j=1}^q x_{i,j} b_{i,j} + \sum_{l=1}^n \sqrt{\sum_{j=1}^q x_{l,j}^2}
\]
where $b_{i,j} \in \{0, 1\}$ for all $i,j$.
\end{itemize}

This list demonstrates that the method is able to deal with differentiable functions,
separable or nonseparable nondifferentiable functions, as well as use several types of atom function in a single problem.

\revisionone{
\subsection{Benchmarking}

In Table~\ref{tab:numexp}, we compare the performance of Algorithm 1 with and without screening, Algorithm~2 with and without screening as well as 2 alternative solvers for 3 problems exhibiting various situations:
\begin{itemize} 
\item Lasso: the nonsmooth function in the Lasso problem is separable;
\item the TV-regularized regression problem has a nonsmooth, nonseparable regularizer whose matrix $A^h$ is sparse;
\item the dual SVM with intercept has a single linear nonseparable constraint.
\end{itemize}

For Algorithm 2, we set the restart with a variable sequence as in \cite{fercoq2018restarting}. We did not tune the algorithmic parameters for each instance.
We evaluate the precision of a primal-dual pair $(x,y)$ as follows. We define the smoothed gap~\cite{tran2018smooth} as
\begin{align*}
\mathcal G_{\beta, \gamma}(x,y,\zeta,\omega) = 
\frac 12 x^\top Q x + f(A^f x) + G(x) + \max_{y'} \langle A^h x, y'\rangle - H^*(y') - \frac{\beta}{2} \norm{y - y'}^2 \\ + H^*(y)   + \frac 12  \omega^\top Q^{\dagger}  \omega + f^*(\zeta)
+ \max_{x'} \langle -(A^h)^\top y-(A^f)^\top \zeta-\omega, x'\rangle - G(x') - \frac{\gamma}{2} \norm{x - x'}^2
\end{align*}
and we choose the positive parameters $\beta$ and $\gamma$ as
\begin{equation*}
\beta = \dist(A^h x, \dom H) \qquad \gamma = \dist(-(A^h)^\top y-(A^f)^\top \zeta-\omega, \dom G^*) \;.
\end{equation*}
It is shown in~\cite{tran2018smooth} that when $\mathcal G_{\beta, \gamma}(x,y,\zeta,\omega)$, $\beta$ and $\gamma$ are small then the objective value and feasibility gaps are small.

For the Lasso problems, we compare our implementations of the coordinate descent method with scikit-learn's coordinate descent \cite{pedregosa2011scikit} and CVXPY's augmented Lagrangian method~\cite{cvxpy} called OSQP. As in Table~\ref{tab:time_comparison}, we have a factor 4 between Alg.\ 1 without screening and scikit-learn. Acceleration and screening allows us to reduce this gap without sacrificing generality.
OSQP is efficient on small problems but is not competitive on larger instances.

For the TV-regularized regression problems, we compare ourself with FISTA where the proximal operator of the total variation is computed inexactly and with LBFGS where the total variation is smoothed with a decreasing smoothing parameter. Those two methods have been implemented for~\cite{dohmatob2014benchmarking}. They manage to solve the problem to an acceptable accuracy in a few hours. As the problem is rather large, we did not run OSQP on it.
For our methods, as $h$ is nonzero, we cannot use variable screening with the current theory. Alg.\ 1 quickly reduces the objective value but fails to get a high precision in a reasonable amount of time. On the other hand, Alg.\ 2 is the quickest among the four solvers tested here.

The third problem we tried is dual support vector machine with intercept.
A very famous solver is libsvm~\cite{chang2011libsvm}, which implements SMO \cite{platt199912}, a 2-coordinate descent method that ensures the feasibility of the constraints at each iteration. The conclusions are similar to what we have seen above. The specialized solver remains a natural choice. OSQP can only solve small instances. Alg.\ 1 has trouble finding high accuracy solutions. Alg.\ 2 is competitive with respect to the specialized solver.

\begin{landscape}
\begin{table}
\small
\revisionone{
\hspace{-2em}\begin{tabular}{|l|l|r|r|r|r|}
\hline
Problem & Dataset & Alternative solver 1 & Alternative solver 2 & Algorithm \ref{algo:main} & Algorithm \ref{algo:accel} \\
\hline
~&&&&&\\[-1ex]
Lasso ($\lambda = 0.1 \norm{X^\top y}_\infty$) & triazines ($J$=186, $I$=60) & Scikit-learn: 0.005s & OSQP: 0.033s & 0.107s; scr: 0.101s & 0.066s; scr: 0.079s\\
$\epsilon=10^{-6}$    & scm1d ($J$=9,803, $I$=280) & Scikit-learn: 24.40s & OSQP: 33.21s & 91.97s; scr: 8.73s & 24.13; scr: 3.21s \\
    & news20.binary ($J$=19,996, $I$=1,355,191) & Scikit-learn: 64.6s & OSQP: $>$2,000s & 267.3s; scr: 114.4s & 169.8s; scr: 130.7s \\[1ex]
\hline
~&&&&&\\[-1ex]
 TV-regularized regression & fMRI \cite{tom2007neural,dohmatob2014benchmarking} $(\alpha_1=\alpha_2 = 5.10^{-3})$&inexact FISTA: 24,341s&  LBFGS+homotopy: 6,893s& $>$25,000s& 2,734s\\
$\epsilon=1$  & ($J$=768, $I$=65,280, $L$=195,840) &&&&
\\[1ex]
\hline
~&&&&&\\[-1ex]
Dual SVM with intercept & ionosphere ($J$=14, $I$=351)& libsvm: 0.04s &OSQP: 0.12s & 3.23s&0.42s \\
$\epsilon=10^{-3}$& leukemia ($J$=7,129, $I$=72) & libsvm: 0.1s &OSQP: 3.2s &25.5s& 0.8s \\
& madelon ($J$=500, $I$=2,600) & libsvm: 50s & OSQP: 37s & 3842s& 170s\\
& gisette ($J$=5,000, $I$=6,000)& libsvm: 70s & OSQP: 936s& 170s & 901s\\
& rcv1 ($J$=47,236, $I$=20,242) & libsvm: 195s & OSQP: Memory error & $>$5000s & 63s \\[1ex]
\hline
\end{tabular}
}
\caption{Time to reach precision $\epsilon$. scr = with gap safe screening. For SVM, we scale the data so that in each column of the data matrix $A$, $A_{i,j} \in [-1,1]$. }
\label{tab:numexp}
\end{table}
\end{landscape}

\section{Conclusion}

This paper introduces a generic coordinate descent solver. The technical challenge behind the implementation is the fundamental need for a compiled language in the low-level operations
that are partial derivative and proximal operator computations. We solved it using pre-defined
atom functions that are combined at run time using a python interface. 

We show how genericity allows us to decouple algorithm development from a particular application problem. As an example, our software can solve at least 12 types of optimization problems on large instances using primal-dual coordinate descent, momentum acceleration, restart and variable screening. 

As future works, apart from keeping the algorithm up to date with the state of the art, we plan to bind our solver with CVXPY in order to simplify further the user experience. 
}

\bibliographystyle{plain}
\bibliography{literature,../../hdr/fercoq}

\end{document}